\RequirePackage{fix-cm}
\documentclass[smallextended]{svjour3-b}       
\smartqed  
\usepackage{graphicx}
\usepackage{multirow}
\newcommand{\D}{\displaystyle}

\begin{document}

\title{GA based robust blind digital watermarking
}


\author{V\'{\i}ctor \'{A}lvarez%
\and Jos\'{e} Andr\'{e}s Armario\and Mar\'{\i}a Dolores Frau\and F\'{e}lix Gudiel\and Mar\'{\i}a Bel\'{e}n G\"{u}emes \and Elena Mart\'{\i}n \and Amparo Osuna
}


\authorrunning{V.Alvarez, J.A.Armario, M.D.Frau, F.Gudiel, M.B.G\"{u}emes, E.Mart\'{\i}n, A.Osuna} 

\institute{V. \'{A}lvarez, J.A. Armario, M.D. Frau, F. Gudiel, E. Mart\'{\i}n, A. Osuna \at
              Dept. Matem\'{a}tica Aplicada 1, ETSII, Avda. Reina Mercedes s/n, 41012, Sevilla, Spain \\
              Tel.: +34-95455-2797, -4386, -4389, -6225, -2798, -2798\\
              Fax: +34-954557878\\
              \email{\{valvarez,armario,mdfrau,gudiel,emartin,aosuna\}@us.es}           
\and M.B. G\"{u}emes \at
              Dept. Algebra, Fac. Matem\'{a}ticas, Avda. Reina Mercedes s/n, 41012, Sevilla, Spain \\
              Tel.: +34-954556969\\
              Fax: +34-954556938\\
              \email{bguemes@us.es}           \and
                         }


\maketitle

\begin{abstract}
A genetic algorithm based robust blind digital watermarking scheme is presented. Starting from a binary image (the original watermark), a genetic algorithm is performed searching for a permutation of this image which is as uncorrelated as possible to the original watermark. The output of the GA is used as our final watermark, so that both security and robustness in the watermarking process is improved. Now, the original cover image is partitioned into non-overlapped square blocks (depending on the size of the watermark image). Then a (possibly extended) Hadamard transform is applied to these blocks, so that one bit information from the watermark image is embedded in each block by modifying the relationship of two coefficients in the transformed matrices. The watermarked image is finally obtained by simply performing the inverse (extended) Hadamard transform on the modified matrices. The experimental results show that our scheme keeps invisibility, security and robustness more likely than other proposals in the literature, thanks to the GA pretreatment.
\end{abstract}

\section{Introduction}

Digital watermarking concerns those methods about how to hide a special mark into digital multimedia data to solve the problems of legal ownership, integrity and authenticity of the original data \cite{CMB02}.

The techniques proposed so far can be grouped into two different approaches, depending on whether the watermark is embedded into the least significant bits (spatial domain approach, \cite{STO94}) or it is embedded attending to the perceptually most significant frequency components of the container image (frequency domain approach, \cite{CKLS97}). Usually one tends to apply techniques of the second type, since spatial domain approaches have relatively low information hiding capacity and, what is more important, can be easily erased by lossy image compression.

Most of frequency domain approaches use discrete wavelet transform (DWT), discrete Fourier transform (DFT) and discrete cosine transform (DCT). Very recently fast Hadamard transform (FHT) has arised as a promising alternative (see \cite{HSTK02} and \cite{ZLZ10} for instance). Interested readers in {\em Hadamard matrices} (square matrices consisting in pairwise orthogonal rows) are referred to \cite{Hor07}.

No matter the processing speed is, watermarking is usually required to muster three conditions: {\em security}, {\em imperceptibility} and {\em robustness}.

Security is concerned with embedding a watermark into a piece of content at an untrusted user device without compromising the security of the watermark key, the watermark or the original (see \cite{LKCV06} for instance).

Perceptibility measures whether perceptible artifacts on the watermarked image are introduced, that is, if the presence of the watermark in the final image is noticeable. This magnitude is measured in terms of the {\em Peak Signal to Noise Ratio}, or $PSNR$ in brief. It is most easily defined via the {\em mean squared error} ($MSE$), so that for images with maximum possible pixel value $range$ (i.e. 255 or 1 depending on whether byte or real storing method is adopted), $PSNR$ is calculated as: \begin{equation}\label{psnr} PSNR=10 \log _{10} \frac{range^2}{MSE}.\end{equation}
Here, for two $m\times n$ monochrome images $K=(k_{i,j})$ and $L=(l_{i,j})$ (where one of the images is considered a noisy approximation of the other), $MSE$ is defined as follows:$$MSE=\frac{1}{mn} \sum _{i=0}^{m-1}\sum _{j=0}^{n-1} (k_{i,j}-l_{i,j})^2.$$


The robustness of a watermark depends on  whether it fails to be detected after unintentional or even malicious transformations (see \cite{VDKP03} for details).
It is usually measured in terms of the {\em Normalized Correlation} ($NC$) between the extracted watermark image $EW=(ew_{i,j})$ (presumably modified) and the original watermark $W=(w_{i,j})$, \begin{equation}\label{nc} NC_{W,EW}= \frac{\D\sum _{i=0}^{m-1} \sum _{j=0}^{n-1} w_{i,j}ew_{i,j}}{\sqrt{\D\sum _{i=0}^{m-1} \sum _{j=0}^{n-1} w_{i,j}^2} \sqrt{\D\sum _{i=0}^{m-1} \sum _{j=0}^{n-1} ew_{i,j}^2}}.\end{equation}


Given a watermarking scheme, it may be straightforwardly improved in terms of security, imperceptibility and robustness by simply introducing some pretreatment to the watermark image in order to destroy space relativity (see \cite{ZXLL09} or \cite{ZLZ10} for instance).

Taking the work in \cite{ZLZ10} as starting point, in this paper  we describe an improved  blind (that is, the original cover image is not needed for extracting the watermark) watermarking scheme, with the following advantages:

\begin{itemize}

\item There is no dependence on the sizes of the watermark and cover images (in \cite{ZLZ10} it is forced to be $1/8$).

\item The trade-off between imperceptibility and robustness is measured in terms of a parameter $b$. The greater $b$ is, the nearer $NC$ is to 1, the smaller $PSNR$ is. Accordingly, the smaller $b$ is, the smaller $NC$ is, the greater $PSNR$ is. Some explanations (beyond simple computational evidence!) will be given in order to justify the optimal value for $b$, depending on the way in which the image is being stored (real or byte representation).

\item Robustness, security and imperceptibility of the scheme are  significantly improved thanks to a pretreatment of the watermark image. We have designed a genetic algorithm (in the sequel, GA in brief), looking for a permutation of the original watermark which is as uncorrelated as possible to it. This GA, equipped with a specific crossover operator specially designed for the occasion, beats usual GA equipped with classical crossover operators concerning permutations problems (such as {\em order 1}, {\em partially mapped} and {\em cycle} crossovers, or {\em edge recombination}).

\end{itemize}

We organize the paper as follows. Section 1 is devoted to introduce the problem of watermarking, and our proposal. The GA looking for permuted images of the watermark with low correlation is described in Section 2. Section 3 is devoted to explain the watermarking scheme. Some executions and examples are showed in Section 4. We include a last section for conclusions and comments.

\section{GA for uncorrelated permuted images}

Given an image $Im$, we want to find a permuted image $PIm$ of $Im$, so that their normalized correlation $NC(Im,PIm)$ is as less as possible.
The key problem here is establishing a method for looking for  permutations of the watermark as uncorrelated as possible. Since one cannot afford to perform an exhaustive search in the full set of permutations, we are dealing with a problem worth of some kind of heuristic solution.

Permutation encoded problems are often used to represent scheduling problems and classic combinatorial optimization problems, such as the Traveling Salesman Problem, Bin Packing or Job Scheduling, to name some. Here, the goal consists in  (or can be solved by) arranging some objects in a certain order.

These optimization problems can seldom be solved by using exhaustive search\-es, due to the large size of the search space. Some kind of heuristic is required instead. For instance, GAs.

In what follows, we assume that the reader is familiar with the general framework of GAs, and their usual elements and characteristics. If necessary, \cite{ES03} (and the references there included) is a good place to get a general overview of the subject.

The GAs which deal with permutations are known as {\em ordering GAs} or simply {\em order-based GAs}.

Although GAs are widely recognized as  powerful and widely applicable optimization methods for string encoded problems,  there is no standard GA for manipulating ordered-list representations. However, several crossover operators and mutations have been suggested in the literature, including {\em Order Crossover} ($OX$), {\em Partially Mapped Crossover} ($PMX$), {\em Cycle Crossover} ($CX$), {\em Edge Recombination} ($ER$), {\em Insert Mutation} ($InsM$), {\em Swap Mutation} ($SwM$), {\em Inversion Mutation} ($InvM$) or {\em Scramble Mutation} $(ScM$), for instance. The interested reader is referred to \cite{ES03} and \cite{PC95} for details and further bibliography.

Here we use a {\em Steady-State} GA, in the sense that three offspring (two coming from crossover, and one more coming from mutation) are generated per generation, which will replace the worst adapted individuals at the moment. This way, we use also {\em elitism}, since we always keep the fittest solution so far. Additionally, we use a ``no duplicates'' policy, in the sense that identical individuals are not allowed to occur in the same generation.

In order to select two individuals for reproduction purposes, we use {\em rank-based} selection, by means of {\em linear ranking}. More concretely, assume that the size of the population is $\mu$. Then sort the population in terms of fitness, so that fittest has rank $\mu$ and worst rank $1$. Now fix a factor $1 \leq s \leq 2$ (we use $s=1.5$ in the sequel). In this circumstances, the probability that the $i^{th}$ individual is selected for reproduction is given by \begin{equation}\label{reproducir} P_i=\frac{2-s}{\mu} + \frac{2(i-1)(s-1)}{\mu (\mu -1)} \end{equation}
Depending on whether the factor $s$ is closer to 1, fitness is accordingly relativized for choosing the individual.

Our population consists of 20 individuals, and every run is limited to 50 generations. The fitness function consists in the normalized correlation (\ref{nc}) between the permuted image and the original watermark, so that the lesser $NC$ is, the fitter an individual is.

Although defining crossover operators for permutation encoding is a difficult task, nevertheless, we have designed a new crossover operator (which we denote simply by $X$), attending to the characteristics of our problem.

Since our watermark images are binary (just black (=0) and white (=1) pixels) $m\times m$ matrices, we can easily encode them as $m^2$-length binary vectors. Furthermore, we can just save the positions in which 1 (analogously, 0) entries are displayed. Assume that the watermark image consists of $k$ white pixels. Then the set of its permuted images is uniquely determined by the set of $k$-subsets of $\{1, \ldots ,m^2\}$.

Given two such different $k$-subsets $S_1$ and $S_2$, we select proportionally positions in $S_1$ and $S_2$ attending to their fitness. Assume that $S_1$ has better fitness than $S_2$. Fix randomly a real number $0.5 \leq r \leq 1$. Then we will get a $\lfloor k\cdot r \rfloor$-subset $S$ of $S_1$, and join this subset with a random $k-\lfloor k\cdot r \rfloor$-subset of $S_2 -S$.

A priori, there is no evidence of which among the crossover operations cited above is more suitable for our purposes.

The comparison in Table \ref{cruces} suggests that our proposed crossover operation $X$ is more successful than traditional ordering crossovers $OX$, $PMX$, $CX$ and $ER$.

Due to the binary character of the watermarks, we have chosen two images with significantly different density of white pixels (the {\em FAMA} logo of the University of Seville and the common {\em Stop} driving signal, see Table \ref{marcas}).

\begin{center}\begin{table}[!ht]\label{marcas}\begin{center}\begin{tabular}{ccc} \includegraphics[width=10mm]{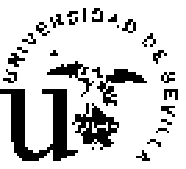} & $\mbox{$\qquad \qquad$}$ &  \includegraphics[width=10mm]{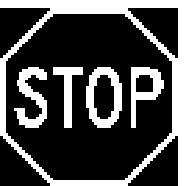} \end{tabular}\end{center} \caption{{\em FAMA} and {\em STOP} watermarks.} \end{table}\end{center}



For each of these images, we have performed 10 runs for each crossover and mutation operators, each of which consisted of 50 generations. The size of every population is fixed in 20 individuals. In the table below, $NC_i$ denotes the best $NC$ value at generation $i$, $Av_i$ denotes the average of $NC_i$ values  along the 10 runs, and $iter$ denotes the average of the generation in which fittest $NC_{50}$ was found along the 10 runs. The optimal values for $NC$ found so far are written in bold. They have been obtained performing a GA consisting of our proposed crossover operator $X$ and the Inverse Mutation operator.


\begin{center}\begin{table*}[!ht]\label{cruces}
$$\begin{array}{c|c|c|c|c|c|c||c|c|c|c|c|}
\multicolumn{2}{c}{}& \multicolumn{5}{c||} {FAMA} & \multicolumn{5}{c|} { STOP}\\ \cline{2-12}
 && NC_0&Av_0&NC_{50}&Av_{50}&iter.& NC_0&Av_0&NC_{50}&Av_{50}&iter.\\   \hline \hline \multirow{4}{*}{OX}&
InsM&0.8048 &0.8076 &0.8021 &0.8037 &36 &0.1872&0.1968&0.1694&0.1803&38\\
\cline{2-12} &  SwM&0.8042&0.8069&0.8009&0.8030&28&0.1794&0.1935&0.1750&0.1819&39\\
\cline{2-12} & InvM&0.8051&0.8077&0.8021&0.8040&35&0.1839&0.1957&0.1727&0.1804&37\\
\cline{2-12}&   ScM&0.8024&0.8068&0.8009&0.8034&29&0.1850&0.1982&0.1694&0.1794&41\\
\hline \hline \multirow{4}{*}{PMX}& InsM&0.8057&0.8078&0.8000&0.8017&43&0.1850&0.1956&0.1627&0.1706&44\\
\cline{2-12} &  SwM&0.8057&0.8075&0.7994&0.8023&43&0.1928&0.1986&0.1694&0.1761&41\\
 \cline{2-12} & InvM&0.8057&0.8071&0.8000&0.8017&44&0.1895&0.1974&0.1594&0.1682&42\\
 \cline{2-12}&   ScM&0.8045&0.8077&0.7976&0.8016&44&0.1672&0.1921&0.1549&0.1668&44\\
 \hline \hline \multirow{4}{*}{CX}& InsM&0.8045&0.8080&0.8036&0.8052&30&0.1884&0.1967&0.1817&0.1877&29\\
\cline{2-12} &  SwM&0.8045&0.8072&0.8042&0.8055&16&0.1906&0.19810&0.1806&0.1898&29\\
\cline{2-12} & InvM&0.8045&0.8073&0.8036&0.8048&23&0.1828&0.1959&0.1750&0.1822&36\\
\cline{2-12}&   ScM&0.8051&0.80786&0.8030&0.8047&21&0.1895&0.1975&0.1761&0.1857&33\\
 \hline \hline \multirow{4}{*}{ER}& InsM&0.8039&0.8074&0.8039&0.8054&30&0.1906&0.1991&0.1806&0.1848&38\\
\cline{2-12} &  SwM&0.8075&0.8082&0.8057&0.8072&21&0.1794&0.1955&0.1794&0.1898&19\\
\cline{2-12} & InvM&0.8054&0.8075&0.8036&0.8053&32&0.1850&0.1936&0.1783&0.1853&16\\
\cline{2-12}&   ScM&0.8063&0.8081&0.8033&0.8055&36&0.1861&0.1966&0.1806&0.1845&33\\
 \hline \hline \multirow{4}{*}{X}& InsM&0.8048&0.8069&0.799&0.8018&43&0.1806&0.1947&0.1560&0.1623&46\\
\cline{2-12} &  SwM&0.8045&0.8081&0.8000&0.8025&40&0.1872&0.1981&0.1505&0.1591&48\\
\cline{2-12} & InvM&0.8057&0.8079&{\bf 0.7952}&0.7993&39&0.1806&0.191&{\bf 0.1159}&0.1343&45\\
\cline{2-12}&   ScM&0.8036&0.8072&0.7994&0.8007&45&0.1884&0.1986&0.1449&0.1581&43\\
 \hline
\end{array}$$
\caption{GA applied to {\em FAMA} and {\em STOP}.}  \end{table*}\end{center}







Notice that there is no interest in timing considerations here, since this is a preprocessing step in our watermarking scheme. Anyway, it is noticeable that every run is performed in just a few seconds.

From Table \ref{cruces} above we conclude that our GA looking for a minimally correlated permuted image is successful, in the sense that normalized correlations of random permuted images (forming the initial populations) are always greater than the local optimum found by the GA, no matter the crossover operator has been selected. Moreover, the crossover $X$ seems to provide fitter permuted images than usual ordering crossovers.

In the following section we describe a watermarking scheme. We claim that pretreatment of the original watermark (so that a minimally correlated permuted image is obtained and used instead), improves the watermarking scheme in an obvious way, not only from the security point of view (no matter one knows the extracting procedure, the extracted watermark will be meaningless), but sometimes (depending on the concrete image, compare Tables 4 
and 6 below) also from the point of view of invisibility, without loose of robustness.


\section{The watermarking scheme}

Grayscale digital images may be stored both in real (the real values of the pixels moving from 0=black to 1=white) or byte (the integer values of the pixels moving from 0=black to 255=white) encoding. Nevertheless, usually byte encoding is preferred, since not every computational system can support working with real encoding (and it requires a discretization step). We will work here with byte encoding, unless stated otherwise.

Let $A$ be the original cover grayscale digital image, encoded as a $n\times n$ matrix with integer entries in $\{0, \ldots, 255\}$.

Let $W$ be the binary watermark, encoded as a $m \times m$ matrix with $0,1$ entries.

\subsection{Watermarking embedding}

The embedding procedure may be detailed as follows:

\begin{itemize}

\item Find a normalized Hadamard matrix $H$ of size $4t$ closest to $\D \lfloor \frac{n}{m}\rfloor$.

\item Partition the original image $A$ into non-overlapped blocks of size $\D \lfloor \frac{n}{m}\rfloor$. Consider the sub-blocks of size $4t \times 4t$ naturally embedded, which we denote by $A_i$, $1 \leq i \leq m^2$.

\item Apply the {\em extended} Hadamard Transform to $A_i$, in order to obtain \begin{equation}\label{had}B_i=\D \frac{H A_i H^T}{4t} \end{equation}

\item Select two entries $b_1$ and $b_2$ in $B_i$ in the same row (or column), say $b_1=B_i(3,3)$ and $b_2=B_i(3,5)$ for instance. Depending on whether the corresponding pixel $i$ in $W$ is 0 or 1, force that $b_2>b_1$ or $b_2<b_1$ accordingly. To this end, fix a value $b$, and take $d=\D \frac{|b_1-b_2|}{2}$. Then set:

\begin{itemize}

\item If $i=0$ and $b_2\leq b_1$ then actualize $b_1^*=b_1-d-b$, $b_2^*=b_2+d+b$.

\item If $i=1$ and $b_2\geq b_1$ then actualize $b_1^*=b_1+d+b$, $b_2^*=b_2-d-b$.

\end{itemize}

\item The watermarked block $A_i^*$ is obtained by the inverse transform of (\ref{had}), \begin{equation} A_i^*=\D \frac{H^T B_i^* H}{4t} \end{equation}

\end{itemize}

At this point, we would like to make two major comments:

\begin{enumerate}

\item Taking a deeper insight in the Hadamard transform (\ref{had}) one deduces that a change $\Delta$ in $b_i$ translates into a change about $\D \lceil \frac{\Delta}{2t} \rceil$ in $A_i^*=\D \frac{H^T B_i^* H}{4t}  $. Hence, in order to get noticeable byte changes, we should take $b \approx t$. In the case of real encoded images, $b$ may be chosen arbitrarily small (at the risk of decreasing the normalized correlation of the extracted watermark).

\item Although the {\em Hadamard Conjecture} about the existence of these matrices in every order $4t$ remains open, there are well known families of Hadamard matrices filling an infinite amount of sizes $4t$ (see \cite{Hor07} for details, and \cite{GPD} for a computer aided generation of Hadamard matrices).

\end{enumerate}

\subsection{Watermarking extraction}

The extraction procedure is just the inverse procedure of embedding.

Let $A_i^*$ be the $i^{th}$-block of the watermarked image. In order to recover the pixel $i$ of the watermark one must simply proceed as follows:

\begin{itemize}

\item By (\ref{had}), form $B_i^*=\D \frac{H A_i^* H^T}{4t}$.

\item Let $b_1$ and $b_2$ be the entries used in the embedding procedure. If $b_2>b_1$ set $i=0$, and $i=1$ otherwise.

\end{itemize}

\section{Experimental results}

We have used 5 different $512 \times 512$ cover images (see Table 
3), and two different $64 \times 64$ watermarks (see Table 
1). We have fixed $H$ to be the $8 \times 8$ Sylvester Hadamard matrix, $$H =\left( \begin{array}{cccccccc} 1& 1& 1& 1& 1& 1& 1& 1\\1& -& 1& -& 1& -& 1& -\\1& 1&-& -& 1& 1& -& -\\1& -& -& 1& 1& -&-& 1\\1& 1& 1&
  1&-&-& -& -\\1& -& 1& -& -& 1& -& 1\\1&
  1& -& -& -& -& 1& 1\\1& -& -& 1& -& 1& 1& -\end{array}\right)$$

\begin{center}\begin{table}[!ht]\label{cover}\begin{center}\begin{tabular}{ccc} \includegraphics[width=20mm]{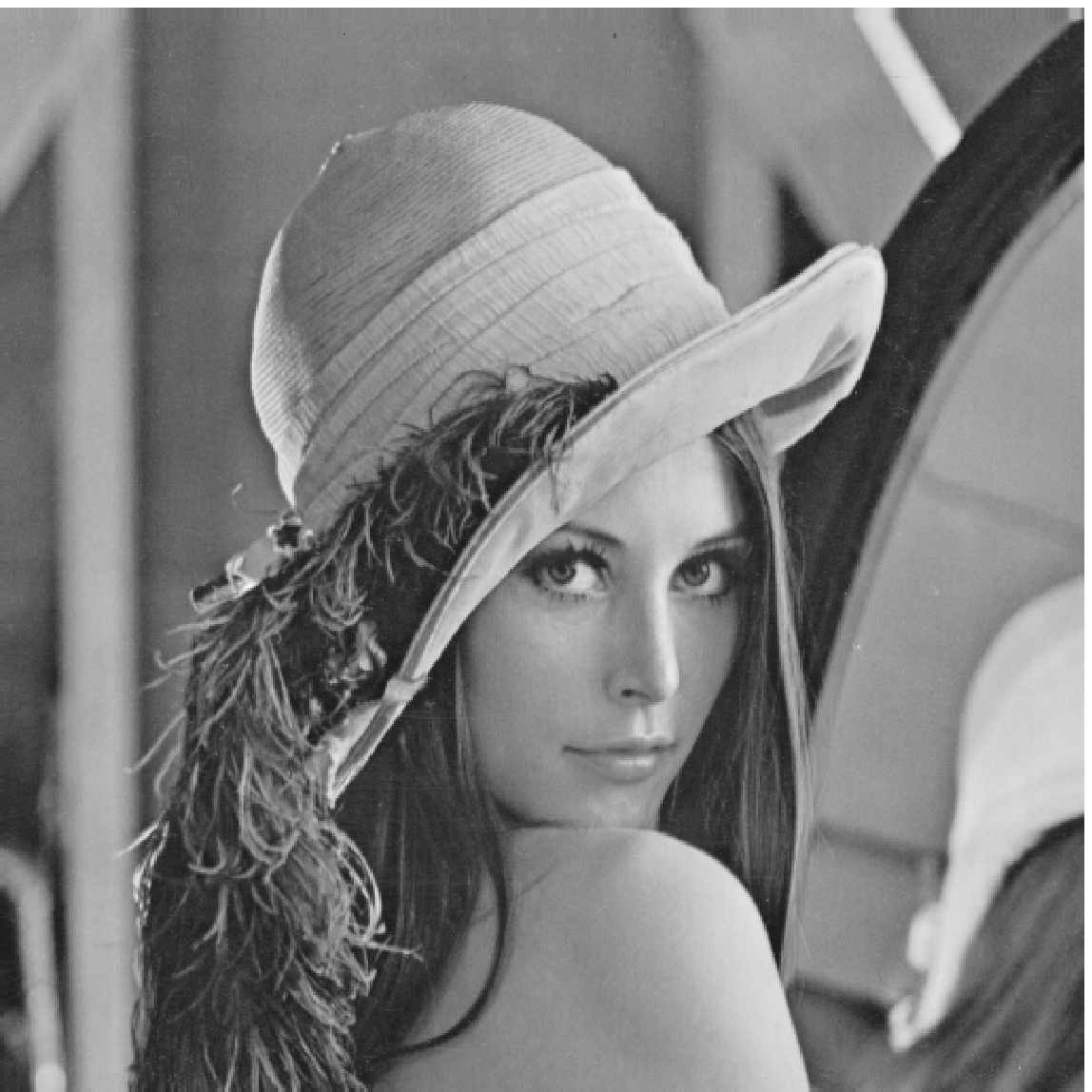} & \includegraphics[width=20mm]{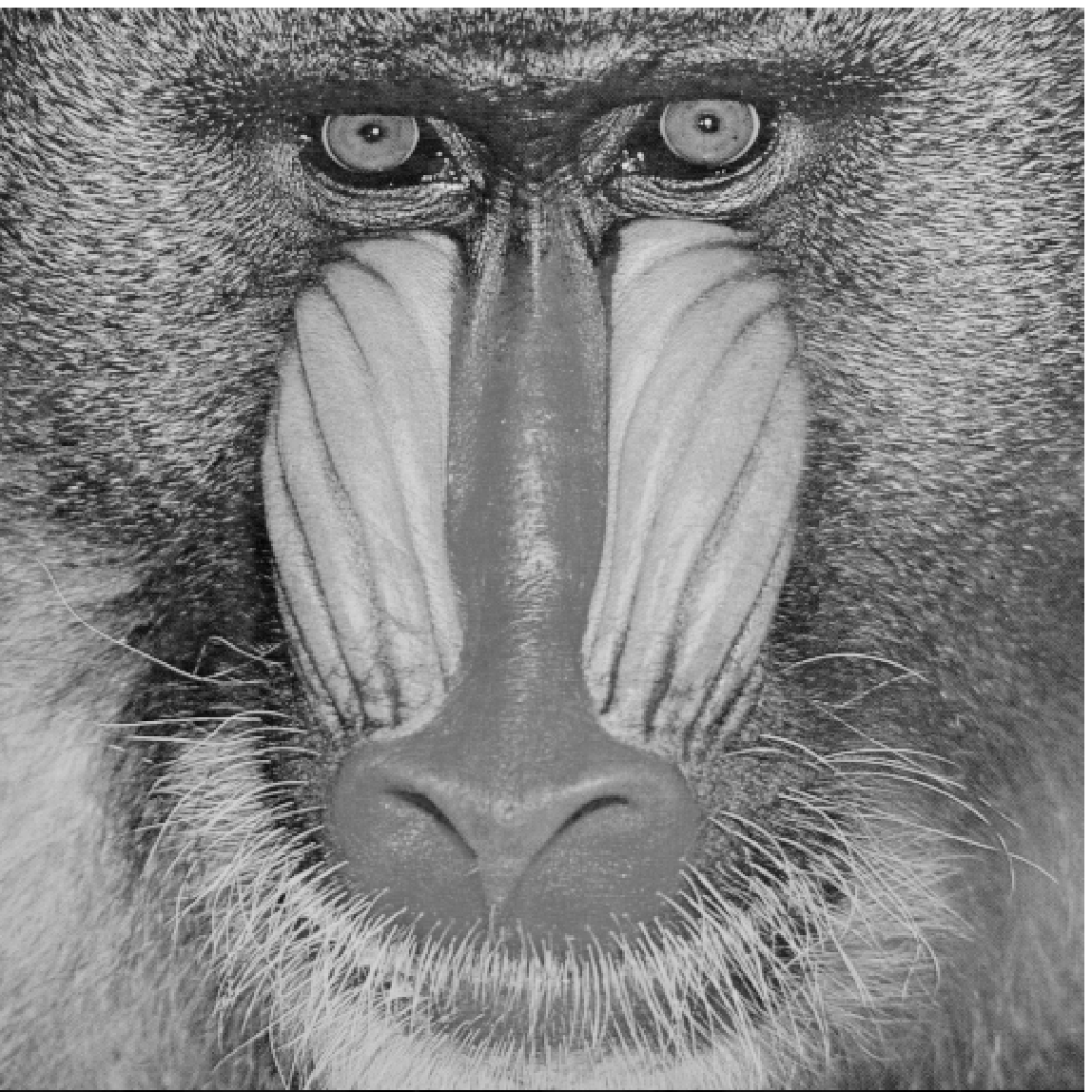} &  \includegraphics[width=20mm]{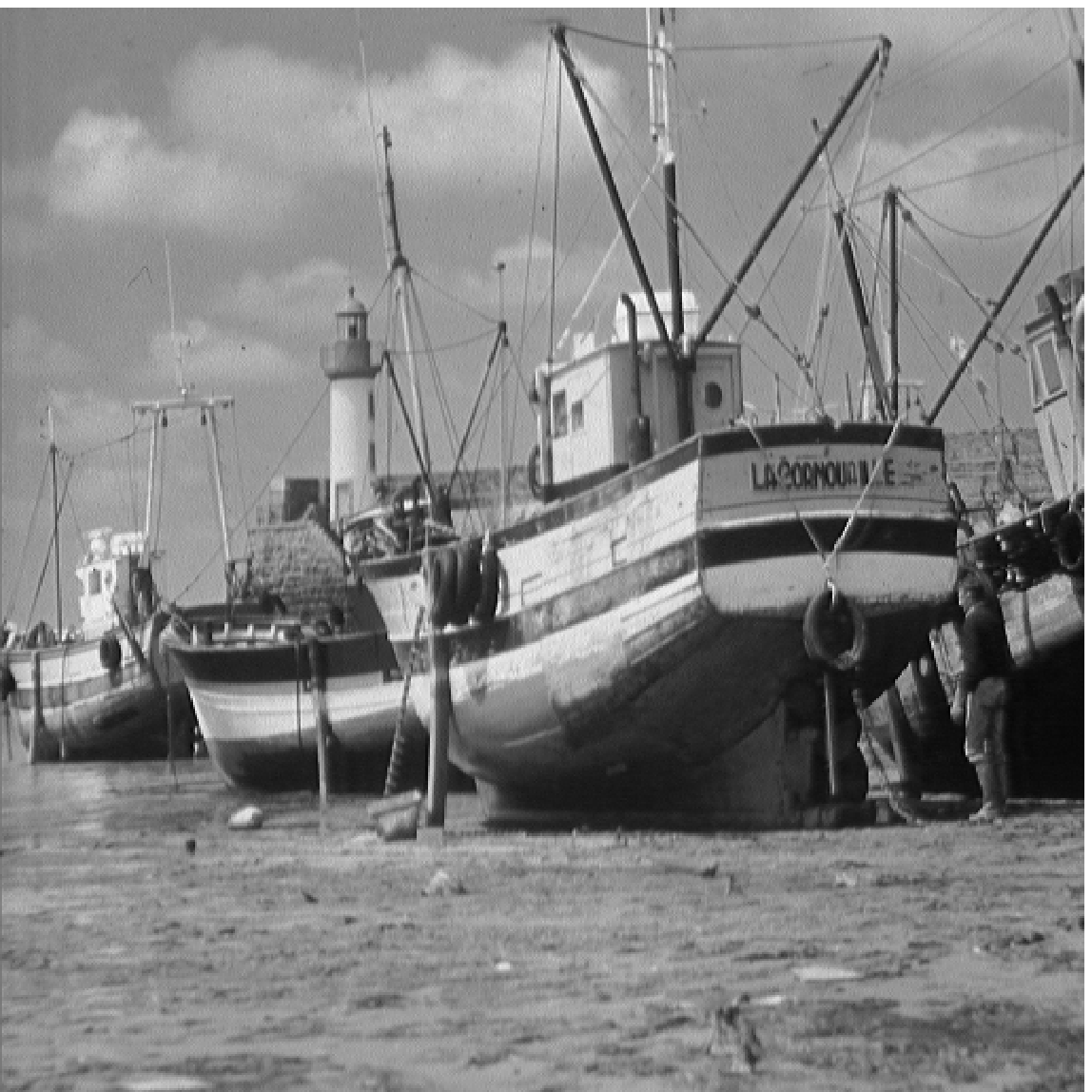}\\ \multicolumn{3}{c}{ \begin{tabular}{cc} \includegraphics[width=20mm]{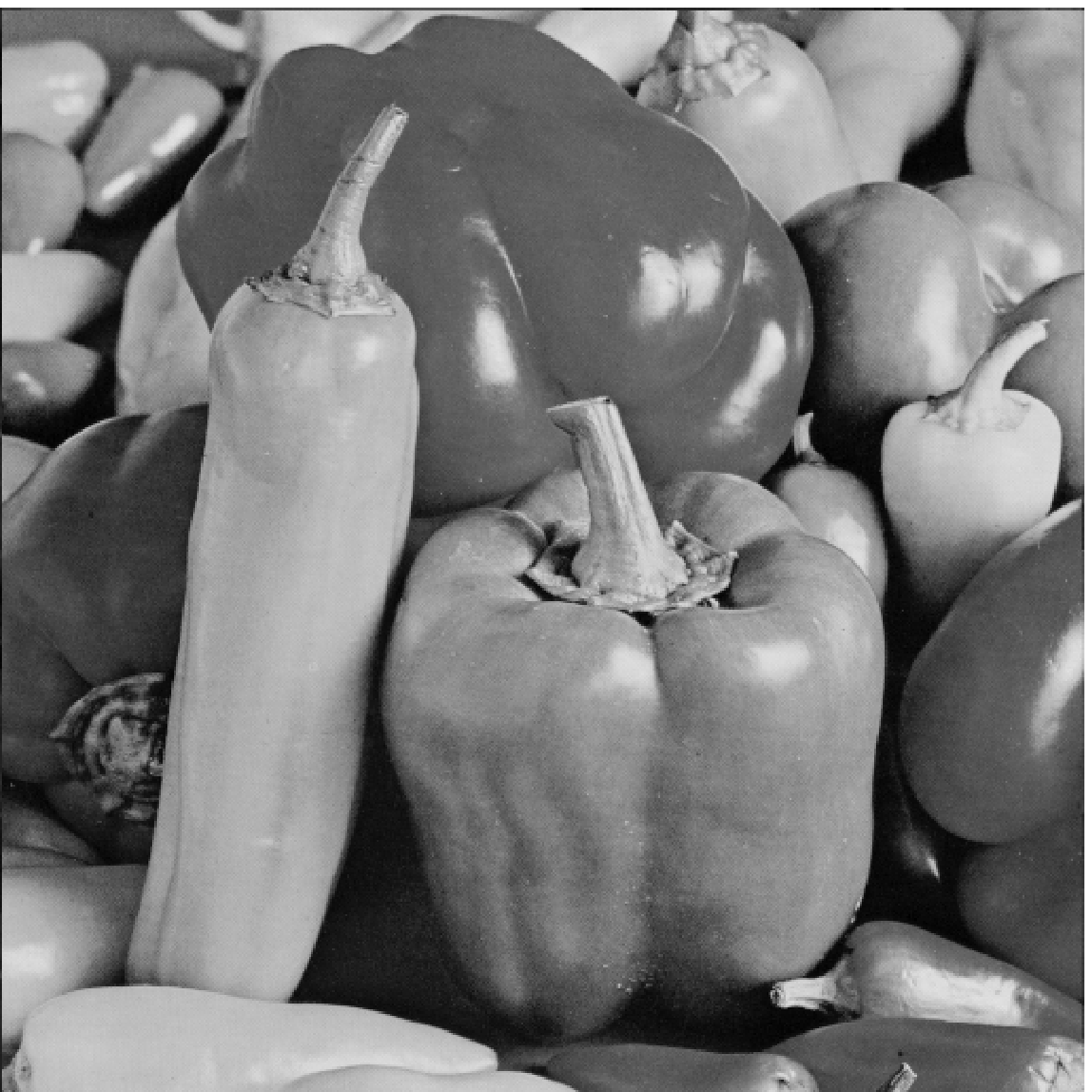}&\includegraphics[width=20mm]{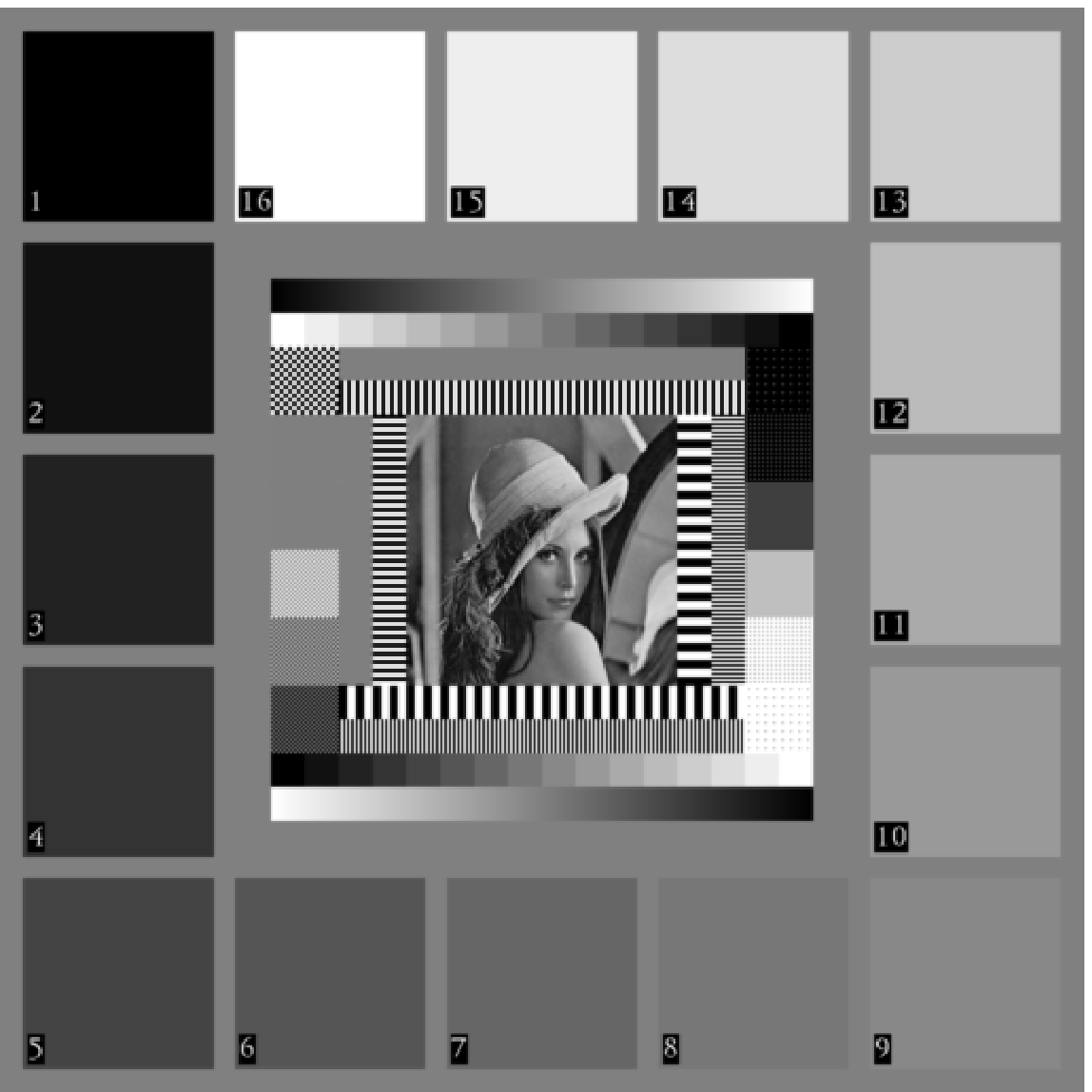} \end{tabular} }\end{tabular}\end{center} \caption{Cover images: Lena, baboon, boats, peppers, testlena.} \end{table}\end{center}

Table 
4 shows the $PSNR$ and $NC$ values obtained for different values of the parameter $b$.

\begin{center}\begin{table*}\label{resultados}
$$\begin{array}{c|c||c|c|c||c|c|c||c|c|c|} \multicolumn{2}{c||}{}& \multicolumn{3}{c||}{Lena}&
 \multicolumn{3}{c||}{Baboon}& \multicolumn{3}{c|}{Boats}\\
\cline{2-11} & b & PSNR& NC&W^*& PSNR& NC&W^*& PSNR& NC&W^*\\
\cline{2-11} \multirow{3}{*}{\includegraphics[width=10mm]{fama.eps}} & 2.01 & 47.1796&1 &\includegraphics[width=10mm]{fama.eps} &40.3194 & 1&\includegraphics[width=10mm]{fama.eps} &46.4154 & 1& \includegraphics[width=10mm]{fama.eps}\\
\cline{2-11}  & 2 &47.2432 & 1&\includegraphics[width=10mm]{fama.eps} &40.34 &1 &\includegraphics[width=10mm]{fama.eps} & 46.4493&1 & \includegraphics[width=10mm]{fama.eps} \\
\cline{2-11}  & 1.99 &47.3031 &0.9976 &\includegraphics[width=10mm]{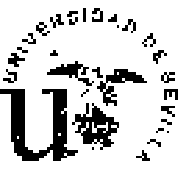} &40.3604 & 0.9989&\includegraphics[width=10mm]{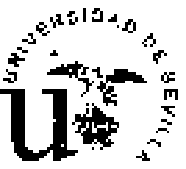} &46.4827 &0.9982 & \includegraphics[width=10mm]{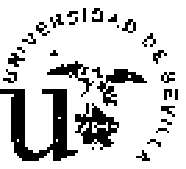} \\
\hline \hline \multirow{3}{*}{\includegraphics[width=10mm]{stop.eps}} & 2.01 & 47.256&1 & \includegraphics[width=10mm]{stop.eps}&40.1841 &1 &\includegraphics[width=10mm]{stop.eps} &46.4054 & 1& \includegraphics[width=10mm]{stop.eps}\\
\cline{2-11}  & 2 &47.2988 &1 &\includegraphics[width=10mm]{stop.eps} &40.2034 &1 &\includegraphics[width=10mm]{stop.eps} & 46.431& 1& \includegraphics[width=10mm]{stop.eps}\\
\cline{2-11}  & 1.99 & 47.3414 & 0.9586&\includegraphics[width=10mm]{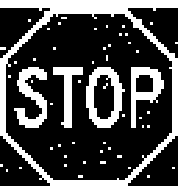} &40.2238 &0.9676 & \includegraphics[width=10mm]{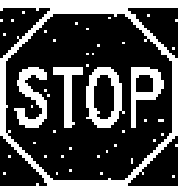}& 46.457& 0.9696& \includegraphics[width=10mm]{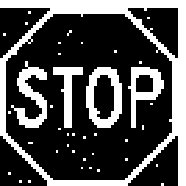} \\ \hline
\end{array}$$\\
%
%
$$\begin{array}{c|c||c|c|c||c|c|c|} \multicolumn{2}{c||}{}& \multicolumn{3}{c||}{Peppers}&
 \multicolumn{3}{c|}{Testlena}\\
\cline{2-8} & b & PSNR& NC&W^*& PSNR& NC&W^*\\
\cline{2-8} \multirow{3}{*}{\includegraphics[width=10mm]{fama.eps}} & 2.01 & 47.256&1 &\includegraphics[width=10mm]{fama.eps}&38.0736 &1 & \includegraphics[width=10mm]{fama.eps}\\
\cline{2-8}  & 2 &48.5626 &1 &\includegraphics[width=10mm]{fama.eps} &38.1818 &0.9711& \includegraphics[width=10mm]{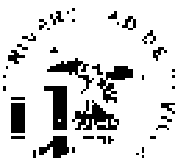} \\
\cline{2-8}  & 1.99 & 48.6202& 0.9982&\includegraphics[width=10mm]{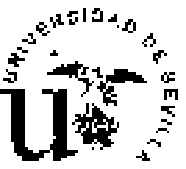} &38.3227 &0.9191 & \includegraphics[width=10mm]{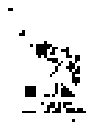}  \\
\hline \hline \multirow{3}{*}{\includegraphics[width=10mm]{stop.eps}} & 2.01 & 48.2392& 1&\includegraphics[width=10mm]{stop.eps} & 42.7214&1 & \includegraphics[width=10mm]{stop.eps} \\
\cline{2-8}  & 2 & 48.3059& 1&\includegraphics[width=10mm]{stop.eps} &43.0388 &0.6448 & \includegraphics[width=10mm]{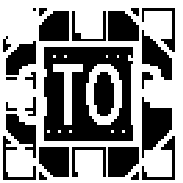} \\
\cline{2-8}  & 1.99 & 48.3734& 0.9480&\includegraphics[width=10mm]{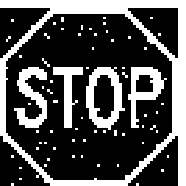} &43.3333 & 0.4993&\includegraphics[width=10mm]{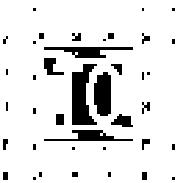}  \\ \hline
\end{array}$$
\caption{$PSNR$ and $NC$ in terms of $b$.} \end{table*}\end{center}





Table 
5 shows that the proposed watermarking scheme is robust under different attacks, such as jpeg compression (with quality factors 80\% and 90\%), Gaussian noise (of mean 0 and variance 0.001) and salt-and-pepper noise (of density 0.01).

\begin{center}\begin{table*}\label{ruidos}
$$\begin{array}{c|c||c|c|c||c|c|c||c|c|c|} \multicolumn{2}{c||}{}& \multicolumn{3}{c||}{Lena}&
 \multicolumn{3}{c||}{Baboon}& \multicolumn{3}{c|}{Boats}\\
\cline{2-11} & Noise & PSNR& NC&W^*& PSNR& NC&W^*& PSNR& NC&W^*\\
\cline{2-11} \multirow{4}{*}{\includegraphics[width=10mm]{fama.eps}} & jpg_{90\%} &39.9389 &0.8758 &\includegraphics[width=10mm]{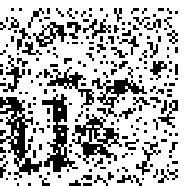} &35.3971 & 0.9057& \includegraphics[width=10mm]{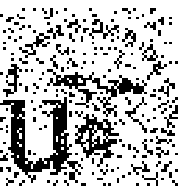}& 38.4282&0.8722 &\includegraphics[width=10mm]{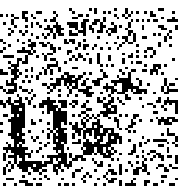} \\
\cline{2-11}  & jpg_{80\%} &38.0837 & 0.7747&\includegraphics[width=10mm]{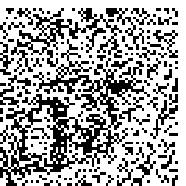} & 31.9305& 0.8313& \includegraphics[width=10mm]{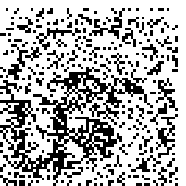}&36.0909 &0.7859 & \includegraphics[width=10mm]{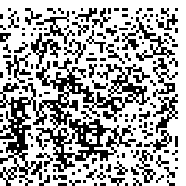} \\
\cline{2-11}  & Gauss. &47.0694 & 0.9969&\includegraphics[width=10mm]{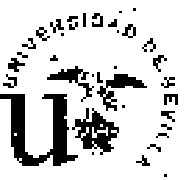} &40.3029 & 0.9978& \includegraphics[width=10mm]{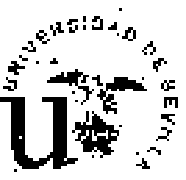}&46.304 &0.9971 &\includegraphics[width=10mm]{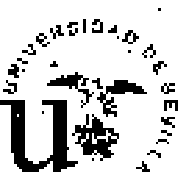}  \\
\cline{2-11}  & S.\& P. &25.4416 & 0.9139&\includegraphics[width=10mm]{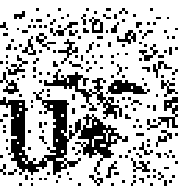} &25.47 &0.9238 &\includegraphics[width=10mm]{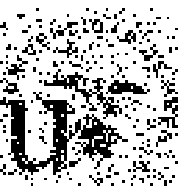} &25.5055 &0.9179 & \includegraphics[width=10mm]{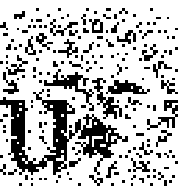}  \\
\hline \hline \multirow{4}{*}{\includegraphics[width=10mm]{stop.eps}} & jpg_{90\%} &39.9704 &0.6477 &\includegraphics[width=10mm]{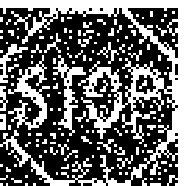} &35.3456 &0.7294 & \includegraphics[width=10mm]{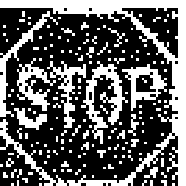}&38.4404 &0.6636 & \includegraphics[width=10mm]{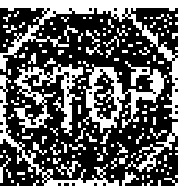} \\
\cline{2-11}  & jpg_{80\%} &38.0843 & 0.4713&\includegraphics[width=10mm]{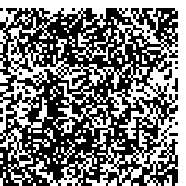} &31.9249 &0.5663 &\includegraphics[width=10mm]{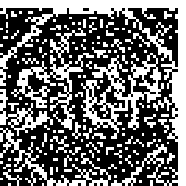} &36.0861 &0.4869 & \includegraphics[width=10mm]{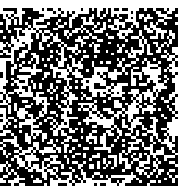} \\
\cline{2-11}  & Gauss. &47.1223 & 0.9729 &\includegraphics[width=10mm]{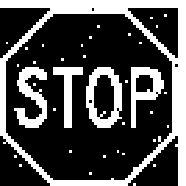} &40.1736 & 0.9884& \includegraphics[width=10mm]{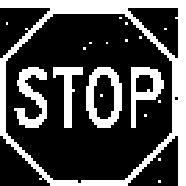}&46.2867 &0.9846 & \includegraphics[width=10mm]{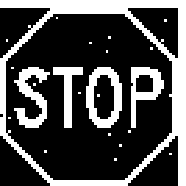} \\
\cline{2-11}  & S.\& P. &25.4419 &0.7591 &\includegraphics[width=10mm]{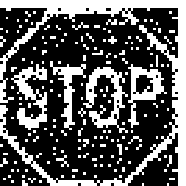} &25.4655 &0.7781 &\includegraphics[width=10mm]{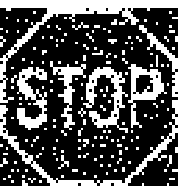} & 25.5054&0.7620 & \includegraphics[width=10mm]{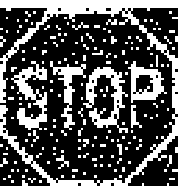} \\ \hline
\end{array}$$
\\
%
%
%
$$\begin{array}{c|c||c|c|c||c|c|c|} \multicolumn{2}{c||}{}& \multicolumn{3}{c||}{Peppers}&
 \multicolumn{3}{c|}{Testlena}\\
\cline{2-8} & Noise & PSNR& NC&W^*& PSNR& NC&W^*\\
\cline{2-8} \multirow{4}{*}{\includegraphics[width=10mm]{fama.eps}} & jpg_{90\%} & 38.4336&0.8739 &\includegraphics[width=10mm]{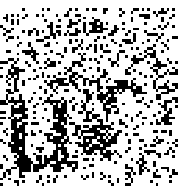} &37.5428 &0.9698 &\includegraphics[width=10mm]{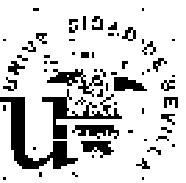} \\
\cline{2-8}  & jpg_{80\%} &36.5783 &0.7764 &\includegraphics[width=10mm]{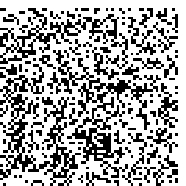} &36.9454 &0.8912 & \includegraphics[width=10mm]{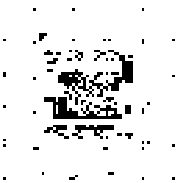} \\
\cline{2-8}  & Gauss. &48.3526 &0.9962 &\includegraphics[width=10mm]{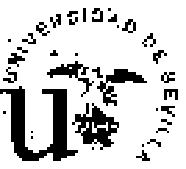} & 38.0697& 0.9941&  \includegraphics[width=10mm]{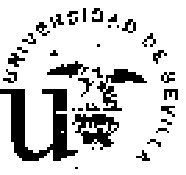} \\
\cline{2-8}  & S.\& P. &25.3601 &0.9140 &\includegraphics[width=10mm]{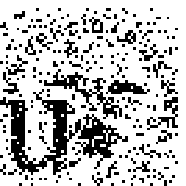} & 24.9004&0.9152 &  \includegraphics[width=10mm]{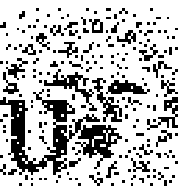} \\
\hline \hline \multirow{4}{*}{\includegraphics[width=10mm]{stop.eps}} & jpg_{90\%} & 38.4118& 0.6660&\includegraphics[width=10mm]{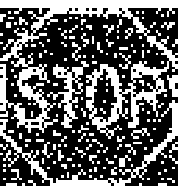} &41.3027 &0.7895 & \includegraphics[width=10mm]{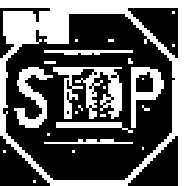}\\
\cline{2-8}  & jpg_{80\%} &36.5486 &0.4955 &\includegraphics[width=10mm]{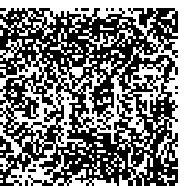} & 39.9323& 0.4621& \includegraphics[width=10mm]{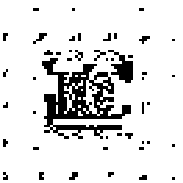} \\
\cline{2-8}  & Gauss. &48.1107 &0.9793 &\includegraphics[width=10mm]{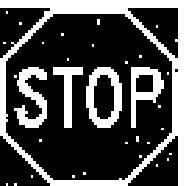} & 42.7086&0.9650 &  \includegraphics[width=10mm]{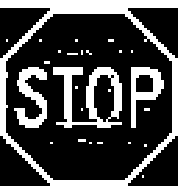} \\
\cline{2-8}  & S.\& P. &25.3589 &0.7548 & \includegraphics[width=10mm]{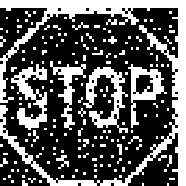}& 25.0392& 0.7564& \includegraphics[width=10mm]{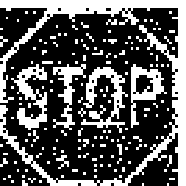}   \\ \hline
\end{array}$$
\caption{$PSNR$ and $NC$ under noises.} \end{table*}\end{center}


\begin{center}\begin{table*}[!ht]\label{nuevopsnr}
$$\begin{array}{c|c||c|c||c|c||c|c||c|c||c|c||} \multicolumn{2}{c}{}& \multicolumn{2}{c||}{Lena}&
 \multicolumn{2}{c||}{Baboon}& \multicolumn{2}{c||}{Boats}& \multicolumn{2}{c||}{Peppers}& \multicolumn{2}{c|}{Testlena}\\
\hline W&PW & PSNR&NC & PSNR&NC & PSNR&NC & PSNR&NC & PSNR&NC\\
\hline\includegraphics[width=10mm]{fama.eps}  & \includegraphics[width=10mm]{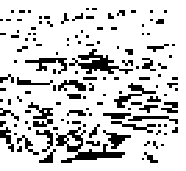}&47.5041&1 &40.2723 &1& 46.4686 &1
 &48.7436 &1 & 43.5181 &1
\\
\hline \hline \includegraphics[width=10mm]{stop.eps}& \includegraphics[width=10mm]{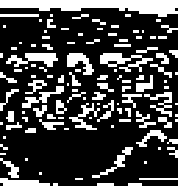}&47.2524 &1&40.4052 &1  & 46.2693&1&
48.1666 &1 & 43.5181&1
\\
\hline
\end{array}$$
\caption{$PSNR$ and $NC$ using GA-watermarks.} \end{table*}\end{center}

Although the proposed watermarking scheme works fine without any need of pretreatment of the watermark, we want to emphasize that permuting the initial watermark ensures that security and imperceptibility (and even robustness to a somewhat lesser degree) are enhanced.

Table 
6 below shows computational evidence of this fact, when we substitute the original watermark $W$ by the permuted images $PW$ provided by the GA described in the previous section (those corresponding to the bold entries in Table 
2).




\section{Conclusions}

In this paper we have described a new blind watermarking scheme. We have shown that pretreatment of the original watermark so that a minimally correlated permuted image is obtained and used instead, improves the watermarking scheme not only from the security point of view, but also from the point of view of imperceptibility, without loose of robustness (see Table 
6).


Although the problem of finding a minimally correlated permuted image from the given watermark is hard, we have designed an order-based Steady-State GA which successfully solve the problem. It includes a proper crossover operator, specifically developed attending to the particular features of our problem. This crossover operator has been shown to beat classical order crossover operators experimentally (see Table 
2).

%
%

\begin{acknowledgements}

All authors are partially supported by FEDER funds via the research projects FQM-016 and P07--FQM--02980 from JJAA and MTM2008-06578 from MICINN (Spain).
\end{acknowledgements}



\end{document}